\documentclass[10pt]{article}
\usepackage{amssymb}
\usepackage{amsmath}
\usepackage{amsbsy}
\newtheorem{lemma}{Lemma}[section]
\newtheorem{proposition}{Proposition}[section]
\newtheorem{theorem}{Theorem}[section]

\def\proclaim#1{\par \bigskip\noindent {\bf #1}\bgroup\it\ }
\def\endproclaim{\egroup\par\bigskip}

\newbox\TempBox \newbox\TempBoxA
\newcommand{\non}{\nonumber \\}
\def\pr{\textsf{P}} 
\def\ep{\textsf{E}} 
\def\Var{\textsf{Var}} 
\def\Cal#1{{\mathcal #1}}

\def\text#1{\mbox{\rm #1}}
\def\overset#1#2{\stackrel{#1}{#2} }

\def\underwiggle 1{
\ifmmode\setbox\TempBox=\hbox{$ 1$}\else\setbox\TempBox=\hbox{
1}\fi \setbox\TempBoxA=\hbox to \wd\TempBox{\hss\char'176\hss}
\rlap{\copy\TempBox}\smash{\lower9pt\hbox{\copy\TempBoxA}} }

\begin{document}




\title{\huge \bf On the rates of  the other law of the logarithm\footnote{Research supported by  National
Natural Science Foundation of China  (No. 10071072)}}
\author{
{\sc   ZHANG Li-Xin\footnote{Department of Mathematics, Zhejiang
University, Hangzhou 310028, China,\newline E-mail:
lxzhang@mail.hz.zj.cn}}
\\
{\em Department of Mathematics, Zhejiang University, Hangzhou
310028, China }}

\date{ }

\maketitle

{\rm

{\sc Abstract.} \quad Let $X$, $X_1$, $X_2$, $\ldots$ be i.i.d.
random variables, and set $S_n=X_1+\ldots + X_n$, $M_n=\max_{k\le
n}|S_k|$, $n\ge 1$. Let $a_n=o(\sqrt{n}/\log n)$. By using the
strong approximation, we prove that: if  $\ep X=0$, $\Var X
=\sigma^2>0$ and $\ep|X|^{2+\epsilon}<\infty$ for some
$\epsilon>0$,  then for any $r>1$,
\begin{eqnarray*}
 \lim_{\epsilon\nearrow 1/\sqrt{r-1}} [\epsilon^{-2}-(r-1)]
\sum_{n=1}^{\infty} n^{r-2} \pr \Big\{M_n \le \epsilon \sigma
\sqrt{\pi^2n/(8\log n)} +a_n\Big\} = \frac 4{\pi}.
\end{eqnarray*}
We also show that the widest $a_n$ is $o(\sqrt{n}/\log n)$.

\bigskip



 \vskip 0.2in

}
\newpage


\section{Introduction and main results.}
\setcounter{equation}{0}

Let $\{X, X_n; n\ge 1\}$ be a sequence of i.i.d random variables
with a common distribution function $F$,  and set
$S_n=\sum_{k=1}^n X_k$, $M_n=\max_{k\le n}|S_k|$, $n\ge 1$. Also
let $\log x=\ln(x\vee e)$, $\log\log x=\log(\log x)$ and
$\phi(x)=\sqrt{\pi^2x/(8\log x)}$. The following is the well known
complete convergence first established by  Hsu and Robbins (1947):
$$ \sum_{n=1}^{\infty} \pr(|S_n|\ge \epsilon n)<\infty, \quad
\epsilon >0
$$
if and only if $\ep X=0$ and $\ep X^2<\infty$. Baum and Katz
(1965) extended this result and proved the following theorem.
\begin{proclaim}{Theorem A} Let $1\le p<2$ and $r\ge p$. Then
$$ \sum_{n=1}^{\infty}n^{r-2} \pr(|S_n|\ge \epsilon n^{1/p})<\infty, \quad
\epsilon >0
$$
if and only if $\ep X=0$ and $\ep |X|^{rp}<\infty$.
\end{proclaim}
Many authors considered  various extensions of the results of
Hsu-Robbins and Baum-Katz. Some of them study the precise
asymptotics of the infinite sums as $\epsilon\to 0$ (c.f. Heyde
(1975), Chen (1978), Sp\u ataru (1999) and Gut and Sp\u ataru
(2000a)). But, this kind of results do not hold for $p=2$.
However, by replacing $n^{1/p}$ by $\sqrt{n\log\log n}$, Gut and
Sp\u ataru (2000b) established an analogous result called the
precise asymptotics of the law of the iterated logarithm. By
replacing $n^{1/p}$ by $\sqrt{n\log n}$, Lai (1974) and Chow and
Lai (1975) consider the following result on the law of the
logarithm.

\begin{proclaim}{Theorem B} Suppose that  $\Var X=\sigma^2$ and
$r\ge 1$. Then the following are equivalent:
$$\sum_{n=1}^{\infty}n^{r-2}\pr(M_n\ge  \epsilon\sqrt{2n\log n})<\infty;
\;\text{ for all } \;\epsilon>\sigma\sqrt{r-1};  $$
$$\sum_{n=1}^{\infty}n^{r-2}\pr(|S_n|\ge   \epsilon\sqrt{2n\log n})<\infty,
\;\text{ for all } \;\epsilon>\sigma\sqrt{r-1}; $$
$$\sum_{n=1}^{\infty}n^{r-2}\pr(|S_n|\ge   \epsilon\sqrt{2n\log n})<\infty,
\;\text{ for some }\; \epsilon>0; $$
$$ \ep X=0\; \text{ and } \; \ep |X|^{2r}/(\log|X|)^r<\infty. $$
\end{proclaim}

Liang, et al (2003)  studied the precise asymptotics of the second
infinite serie in Theorem B for $1<r<3/2$ under the condition
$\ep|X|^{2r}<\infty$. Zhang (2003) studied all the cases of $r>1$
and obtained the sufficient and necessary condition for such kind
of results to hold.

By a small deviation theorem of Mogul'ski{\sc \u i} (1974) (c.f.,
Lemma \ref{lem1}), it is easy to get the following results on the
other law of the logarithm.
\begin{theorem}\label{th0} Suppose that  $\ep X=0$, $\Var X=\sigma^2$ and
$r> 1$. Then
$$\sum_{n=1}^{\infty}n^{r-2}(\log n)^a \pr(M_n\le  \epsilon \sqrt{n/\log n})<\infty
\;\text{ for all } \;\epsilon<\sigma\sqrt{\frac{\pi^2}{8(r-1)}}
$$
and
$$\sum_{n=1}^{\infty}n^{r-2}(\log n)^a\pr(M_n\le  \epsilon \sqrt{n/\log n})=\infty
\;\text{ for all } \;\epsilon>\sigma\sqrt{\frac{\pi^2}{8(r-1)}}.
$$
\end{theorem}

The purpose of this paper is to consider the precise asymptotics
of the infinite series in Theorem \ref{th0} for all $r>1$.  Here
is our main result.

\begin{theorem} \label{th1}
Let $r>1$ and $a>-1$ and let $a_n(\epsilon)$ be a function of
$\epsilon$ such that
\begin{eqnarray} \label{co1.1}
a_n(\epsilon) \log n \to \tau \; \text{ as } \; n\to \infty \text{
and } \epsilon \nearrow 1/\sqrt{r-1}.
\end{eqnarray}
Assume that
\begin{equation}\label{eqT1.0}
\ep X=0, \quad \ep X^2=\sigma^2 \; (0<\sigma<\infty) \; \text{ and
} \; \ep[|X|^{2+\epsilon}]<\infty, \quad \text{ for some }
0<\epsilon<1.
\end{equation}
 Then
\begin{eqnarray}\label{eqT1.1}
& & \lim_{\epsilon\nearrow 1/\sqrt{r-1}}
[\epsilon^{-2}-(r-1)]^{a+1} \sum_{n=1}^{\infty} n^{r-2}(\log n)^a
\pr \Big\{M_n \le \sigma \phi(n) (\epsilon +a_n(\epsilon))\Big\}
\non & & \qquad \qquad \qquad=\frac {4}{\pi }\exp\{2\tau
(r-1)^{3/2}\}\Gamma(a+1).
\end{eqnarray}
Here, $\Gamma(\cdot)$ is a gamma function. Conversely, if
(\ref{eqT1.1}) holds for some $r>1$, $a>-1$ and $\epsilon>0$, then
$\ep X=0$ and $\Var X=\sigma^2$.
\end{theorem}

Also, we have a refinement of Theorem \ref{th0} as
\begin{theorem} \label{th3}
Let $r>1$ and $a$ be two real numbers.  Suppose that the condition
(\ref{eqT1.1}) is satisfied, then for any eventually
non-decreasing $\psi:[1,\infty)\to (0, \infty)$,
\begin{eqnarray} \label{eqT3.1}
& & \sum_{n=1}^{\infty} n^{r-2}(\log n)^a \pr \Big\{ M_n\le
\sigma\sqrt{ \pi^2 n/(8 \psi(n))} \Big\}<\infty \text{ or }
=\infty \non
 & & \text{ according as }
\sum_{n=1}^{\infty} n^{r-2}(\log n)^a \exp\{-\psi(n)\}<\infty
\text{ or } = \infty.
\end{eqnarray}
\end{theorem}

 We conjecture
that (\ref{eqT1.1}) is true whenever $\ep X=0$, $\Var
X=\sigma^2>0$ and $\ep X^2I\{|X|\ge t\}=o\big((\log t)^{-1}\big)$
as $t\to \infty$, and (\ref{eqT3.1}) is true whenever $\ep X=0$,
$\Var X=\sigma^2>0$ and $\ep X^2I\{|X|\ge t\}=O\big((\log
t)^{-1}\big)$ as $t\to \infty$.

The proofs of Theorems \ref{th1} and \ref{th3} are given in
Section 4. Before that, we first verify (\ref{eqT1.1}) under the
assumption that $F$ is the normal distribution in Section 2, after
which, by using the strong approximation method, we then show that
the probability in (\ref{eqT1.1})
 can be replaced by those for
normal random variables in Section 3. Throughout this paper, we
let $K(\alpha,\beta,\cdots)$, $C(\alpha,\beta,\cdots)$ etc denote
positive constants which depend on $\alpha,\beta, \cdots$ only,
whose values can differ in different places. The notation $a_n\sim
b_n$ means that $a_n/b_n\to 1$, and $a_n\approx b_n$ means that
$C^{-1} b_n\le a_n \le C b_n$ for some $c>0$ and all $n$ large
enough.


\section{Normal cases.}
\setcounter{equation}{0}

In this section, we prove Theorem \ref{th1} in the case that $\{X,
X_n; n\ge 1\}$ are normal random variables. Let $\{W(t); t\ge 0\}$
be a standard Wiener process and $N$ a standard normal variable.
Our result is as follows.

\begin{proposition}\label{prop2.1}
Let $r>1$ and $a>-1$ and let $a_n(\epsilon)$ be a function of
$\epsilon$ satisfying (\ref{co1.1}).  Then
\begin{eqnarray}\label{eqprop2.1.1}
& & \lim_{\epsilon\nearrow 1/\sqrt{r-1}}
[\epsilon^{-2}-(r-1)]^{a+1} \sum_{n=1}^{\infty}n^{r-2} (\log n)^a
\pr \Big\{\sup_{0\le s\le 1}|W(s)| \le \sqrt{\pi^2/(8\log n)}
(\epsilon +a_n(\epsilon))\Big\} \non & & \qquad \qquad
\qquad=\frac {4}{\pi }\exp\{2\tau (r-1)^{3/2}\}\Gamma(a+1).
\end{eqnarray}
\end{proposition}

The following lemma will be used in the proofs.

\begin{lemma}\label{lem2.1}
Let $\{W(t); t\ge 0\}$ be a standard Wiener process. Then for all $x>0$,
\begin{eqnarray*}\label{eqL2.1.1}
 \pr \big(\sup_{0\le s\le 1}|W(s)| \le x\big) = \frac{4}{\pi}\sum_{k=0}^{\infty}
 \frac{(-1)^k}{2k+1}\exp\big\{-\frac{\pi^2(2k+1)^2}{8 x^2}\big\}.
\end{eqnarray*}
In particular,
\begin{eqnarray}\label{eqL2.1.1}
\frac{2}{\pi}
 \exp\big\{-\frac{\pi^2}{8 x^2}\big\}\le \pr \big(\sup_{0\le s\le 1}|W(s)| \le x\big) \le \frac{4}{\pi}
 \exp\big\{-\frac{\pi^2}{8 x^2}\big\}
\end{eqnarray}
and
$$
\pr\big( \sup_{0\le s\le 1}|W(s)|\le x \big) \sim
\frac{4}{\pi}\exp\big\{-\frac{\pi^2}{8 x^2}\big\}\; \text{ as } \;
x\to 0.
$$
\end{lemma}
{\bf Proof.} It is well known. See Ciesielaki and Taylor (1962).

\begin{lemma}\label{lem2.2} Let $\alpha_n(\epsilon)>0$,
$\beta_n(\epsilon)>0$  and $f(\epsilon)>0$ satisfying
$$ \alpha_n(\epsilon)\sim \beta_n(\epsilon) \quad \text{ as }\;
n\to \infty \; \text{ and }\; \epsilon\to\epsilon_0, $$
$$   f(\epsilon)\beta_n(\epsilon)\to 0 \;
\text{ as } \epsilon\to\epsilon_0, \forall n. $$
Then
$$ \limsup_{\epsilon\to \epsilon_0}(\liminf_{\epsilon\to
\epsilon_0})f(\epsilon)\sum_{n=1}^{\infty}\alpha_n(\epsilon)=
 \limsup_{\epsilon\to \epsilon_0}(\liminf_{\epsilon\to
\epsilon_0})f(\epsilon)\sum_{n=1}^{\infty}\beta_n(\epsilon).
$$
\end{lemma}
{\bf Proof.} Easy.

\bigskip

Now, we turn to prove the proposition \ref{prop2.1}.
 By Lemma
\ref{lem2.1} and the Condition (\ref{co1.1}) we have
\begin{eqnarray*}
& &\pr \Big\{\sup_{0\le s\le 1}|W(s)|
   \ge \sqrt{\frac{\pi^2}{8\log n}}
(\epsilon +a_n(\epsilon))\Big\} \sim \frac{4}{\pi}\Big\{
-\frac{\log n}{(\epsilon +a_n(\epsilon))^2}\Big\}
\\
& & \qquad =  \frac{4}{\pi}\Big\{ -\frac{\log n}{\epsilon^2 +2
\epsilon a_n(\epsilon)+a_n(\epsilon)^2}\Big\}   \sim
\frac{4}{\pi}\Big\{ -\frac{\log n}{\epsilon^2}\Big\}
\exp\Big\{\frac{2}{\epsilon^3}  a_n(\epsilon)\log n\Big\}  \\
& & \qquad \sim  \frac{4}{\pi}\Big\{ -\frac{\log
n}{\epsilon^2}\Big\} \exp\Big\{2\tau (r-1)^{3/2} \Big\}
\end{eqnarray*}
as $n\to \infty$ and $\epsilon\nearrow 1/\sqrt{r-1}$. We conclude
that
\begin{eqnarray*}
& & \lim_{\epsilon\nearrow1/\sqrt{r-1}}[\epsilon^{-2}-(r-1)]^{a+1}
\sum_{n=1}^{\infty} n^{r-2} (\log n)^a \pr \Big\{\sup_{0\le s\le
1}|W(s)|
   \ge \sqrt{\pi^2/(8\log n)}
(\epsilon +a_n(\epsilon))\Big\}
\\
&=&\lim_{\epsilon\nearrow1/\sqrt{r-1}}[\epsilon^{-2}-(r-1)]^{a+1}
\sum_{n=1}^{\infty} n^{r-2} (\log n)^a \exp\big\{-\frac{\log
n}{\epsilon^2} \big\} \frac {4}{\pi }\exp\{2\tau (r-1)^{3/2} \}
\\
&& \qquad\qquad (\text{ by Lemma \ref{lem2.2}})
\\
&=&\lim_{\epsilon\nearrow1/\sqrt{r-1}}[\epsilon^{-2}-(r-1)]^{a+1}
\sum_{n=1}^{\infty} \int_n^{n+1}x^{r-2}(\log x)^a
\exp\big\{-\frac{\log x}{\epsilon^2} \big\}dx \frac {4}{\pi
}\exp\{2\tau (r-1)^{3/2} \}
\\
&& \qquad\qquad (\text{ by Lemma \ref{lem2.2}})
\\
&=&\lim_{\epsilon\nearrow1/\sqrt{r-1}}[\epsilon^{-2}-(r-1)]^{a+1}
 \int_e^{\infty}x^{r-2}(\log x)^a\exp\big\{-\frac{\log x}{\epsilon^2} \big\}dx \frac {4}{\pi
}\exp\{2\tau (r-1)^{3/2} \}
\\
&=&\lim_{\epsilon\nearrow1/\sqrt{r-1}}[\epsilon^{-2}-(r-1)]^{a+1}
 \int_1^{\infty}y^a
 \exp\Big\{-[\epsilon^{-2}-(r-1)]y\Big\}dy\frac {4}{\pi
}\exp\{2\tau (r-1)^{3/2} \}
\\
&=&\lim_{\epsilon\nearrow1/\sqrt{r-1}}
 \int_{\epsilon^{-2}-(r-1)}^{\infty}y^a e^y\frac {4}{\pi
}\exp\{2\tau (r-1)^{3/2} \}=\frac {4}{\pi }\exp\{2\tau (r-1)^{3/2}
\}\Gamma(a+1) .
\end{eqnarray*}
(\ref{eqprop2.1.1}) is proved.


\section{  Approximation.}
\setcounter{equation}{0}

The purpose of this section is to use strong approximation and
Feller's (1945)
 and Einmahl's (1989)
truncation methods to show that the probability in (\ref{eqT1.1})
for $M_n$ can be approximated by those for $\sqrt{n}\sup_{0\le
s\le 1}|W(s)|$.

Suppose that $\ep X=0$ and $\ep X^2=\sigma^2<\infty$. Without
losing of generality, we assume that $\sigma=1$.
 Let $0<p<1/2$. For each $n$ and $1\le j\le n$,
we let $X_{nj}^{\prime}=X_{nj}I\{|X_j|\le \sqrt{n}/n^p \}$,
$X_{nj}^{\ast}=X_{nj}^{\prime}-\ep [X_{nj}^{\prime}]$,
$S_{nj}^{\prime}=\sum_{i=1}^jX_{nj}^{\prime}$,
$S_{nj}^{\ast}=\sum_{i=1}^jX_{nj}^{\ast}$, $M_n^{\ast}=\max_{k\le
n}|S_{nk}^{\ast}|$ and $B_n=\sum_{k=1}^n \Var(X_{nk}^{\ast})$. The
following  proposition is the main result of this section.

\begin{proposition}\label{prop3.1}
Suppose $\ep[|X|^{2+\epsilon}]<\infty$ for some $0<\epsilon<1$.
Let $a>-1$, $r>1$ and $0<p<\epsilon/(4(2+\epsilon))$. Then there
exist $\delta>0$ and a sequence of positive numbers $\{q_n \}$
such that
\begin{eqnarray} \label{eqprop3.1.1}
&&\pr\Big\{\sup_{0\le s\le 1}|W(s)| \le \epsilon \sqrt{\frac{\pi^2
}{8\log n} }
 + \frac{5}{(\log n)^2} \Big\}
-q_n \non &\le& \pr\Big\{M_n \le \epsilon \phi(n) \Big\}
\non
&\le&\pr\Big\{\sup_{0\le s\le 1}|W(s)| \le \epsilon
\sqrt{\frac{\pi^2 }{8\log n} }
 -
\frac{5}{(\log n)^2} \Big\} +q_n, \non & &  \qquad \forall
\epsilon\in (\frac 1{\sqrt{r-1}}-\delta, \frac
1{\sqrt{r-1}}+\delta), \quad n\ge 1
\end{eqnarray}
and
\begin{eqnarray}\label{eqprop3.1.2}
\sum_{n=1}^{\infty}n^{r-2} (\log n)^a  q_n \le
K(r,a,p,\epsilon,\delta)<\infty.
\end{eqnarray}
\end{proposition}

To show this result, we need some lemmas.

\begin{lemma}\label{lem1}
For any $x>0$ and $0<\delta<1$, there exists a positive constant
$C=C(x,\delta)$ such that
\begin{itemize}
\item[{(a)}] $C^{-1}\exp\{-\frac{1+\delta}{x^2}\log n\}
 \le \pr\{M_n^{\ast}\le x\phi(n)\}
 \le C\exp\{-\frac{1-\delta}{x^2}\log  n\}$,
\item[{(b)}] $C^{-1}\exp\{-\frac{1+\delta}{x^2}\log  n\}
 \le \pr\{M_n \le x\phi(n)\}
 \le C\exp\{-\frac{1-\delta}{x^2}\log  n\}$.
\end{itemize}
\end{lemma}

{\bf Proof.}  This lemma is so-called small deviation theorem. It
follows from Theorem 2 of Shao (1995) by noting that $B_n\sim n$.
(See also Shao 1991).

\begin{lemma}\label{lem2}
For any $x>0$, $A>0$ and $0<\delta<1$, there exists a  positive
constant $C=C(x,\delta)$ such that
$$ \pr\{\max_{k\le q}|S_k+z_1|\vee
\max_{q<k\le n}|S_k+z_2|\le x\phi(n)\} \le
C\exp\big\{-\frac{1-\delta}{x^2}\log  n\big\} $$ holds uniformly
in $|z_1|\le A\phi(n)$, $|z_2|\le A\phi(n)$ and $1\le q\le n$.
\end{lemma}

{\bf Proof.} Without losing of generality, we can assume that
$0<\delta<\frac 1{8^8 5!}$. We follow the lines of the proof of
(17) in Shao (1995). Write $x_n=x\phi(n)$. Put $M=\delta^{-2}$.
For fixed $n$, define $m_0=0$,
$$ m_i=\max\{j: j\le Mi x_n^2\},
 \quad \text{ for } i\le l:=\max\{i: m_i\le n-1\} $$
and $m_{l+1}=n$. It is easily seen that
$$(1-\delta/4)M x_n^2\le m_i-m_{i-1}\le (1+\delta/4) M x_n^2 $$
and
$$\frac{n}{M x_n^2}-1\le l \le \frac{n}{M x_n^2} $$
provided $n$ is sufficiently large. From Lemmas 3 and 1 of Shao
(1995) and the Anderson's inequality, it follows that, there
exists an integer $n_0$ such that $\forall n\ge n_0$, $\forall
1\le j\le l$, $\forall |y|\le x_n+A\phi(n)$, $\forall |y_j|\le
A\phi(n)$,
\begin{eqnarray*}
& & \pr\Big(\max_{k\le m_j-m_{j-1}}|S_{m_{j-1}+k}-S_{m_{j-1}}+y_j+y|\le x_n\Big) \\
&\le & e^{-3M}
 +\pr\Big(\sup_{0\le s\le 1}|W(s)+(y_j+y)(m_j-m_{j-1})^{-1/2}|\le x_n(m_j-m_{j-1})^{-1/2}\Big) \\
&\le & e^{-3M}
 +\pr\Big(\sup_{0\le s\le 1}|W(s)|\le x_n(m_j-m_{j-1})^{-1/2}\Big) \\
&\le & e^{-3M} +4\exp\Big(-\frac{\pi^2(m_j-m_{j-1})}{8x_n^2}\Big)
\le e^{-3M}
+4\exp\Big(-\frac{\pi^2 M(1-\delta/4)}{8}\Big) \\
&\le &\frac 12e^{-2M} +\frac 12 \exp\Big(-\frac{\pi^2
M(1-\delta/2)}{8}\Big) \le \exp\Big(-\frac{\pi^2
M(1-\delta/2)}{8}\Big),
\end{eqnarray*}
where $\{W(t); t\ge 0\}$ is a standard Wiener process. Obviously,
there exists an $i$ such that $m_{i-1}<q\le m_i$. Let $y_j=z_1$ if
$j\le i-1$, and $z_2$ if $j\ge i$. Then
\begin{eqnarray*}
& & \pr\big( \max_{k\le q}|S_k+z_1|\vee
  \max_{q<k\le n}|S_k+z_2|\le x\phi(n)\big)
\le \pr\big( \max_{j\le l,j\ne i}\max_{m_{j-1}<k\le m_j}|S_k+y_j|
 \le x_n \big) \\
&= &\ep\Big\{I\{ \max_{j\le l-1,j\ne i}\max_{m_{j-1}<k\le
m_j}|S_k+y_j|
 \le x_n\}I\{ \max_{m_{l-1}<k\le m_l, l\ne i}|S_k+y_j|
 \le x_n\} \Big\} \\
&\le &\ep\Big\{I\{ \max_{j\le l-1,j\ne i}\max_{m_{j-1}<k\le
m_j}|S_k+y_j|
 \le x_n\}\\
& & \qquad \times \ep\big[I\{ \max_{m_{l-1}<k\le m_l}|S_k+y_j|
 \le x_n\}|S_k, k\le m_{l-1}\big] \Big\} \\
&= &\int_{-x_n-y_{l-1}}^{x_n-y_{l-1}}
 \pr\Big( \max_{m_{l-1}<k\le m_l}|S_k-S_{m_{l-1}}+y+y_l|
 \le x_n\Big) \\
& & \qquad \quad d \pr\Big( \max_{j\le l-1,j\ne i}\max_{m_{j-1}<k\le m_j}|S_k+y_j| \le x_n, S_{m_{l-1}}<y\Big) \\
&\le& \exp\Big(-\frac{\pi^2 M(1-\delta/2)}{8}\Big) \pr\Big(
\max_{j\le l-1,j\ne i}\max_{m_{j-1}<k\le m_j}|S_k+y_j|
   \le x_n\Big) \\
&\le &\ldots \le \exp\Big(-\frac{\pi^2 M(1-\delta/2)(l-1)}{8}\Big) \\
&\le & C \exp\Big(-\frac{(1-\delta)\pi^2 n}{8x_n^2}\Big) \le C
\exp\Big(-\frac{(1-\delta)}{x^2}\log  n\Big).
\end{eqnarray*}
For $n\le n_0$, it is obvious that
$$\pr\big( \max_{k\le q}|S_k+z_1|\vee
  \max_{k<k\le n}|S_k+z_2|\le x\phi(n)\big) \le 1
\le C \exp\Big(-\frac{(1-\delta)}{x^2}\log  n\Big). $$ Lemma 2 is
proved.

\begin{lemma}\label{lem3}
Define $\Delta_n=\max_{k\le n}|S_{nk}^{\ast}-S_k|$.  Suppose
$\ep[|X|^{2+\epsilon}]<\infty$ for some $0<\epsilon<1$. Let
$a>-1$, $r>1$ and $0<p<\epsilon/(4(2+\epsilon))$. Then there exist
constants $\delta_0=\delta_0(r,a,p,\epsilon)>0$ and
$K=K(r,a,p,\epsilon,\delta_0)$ such that $\forall
0<\delta<\delta_0$,
$$\sum_{n=1}^{\infty} n^{r-2} (\log n)^a  I_n\le K <\infty, $$
where
$$I_n=\pr\Big(\Delta_n\ge \sqrt{n}/(\log  n)^2,
 M_n^{\ast}\le \phi(n)\big(\frac 1{\sqrt{r-1}}+\delta\big)\Big). $$
\end{lemma}
{\bf Proof.} It is sufficient to show that
\begin{eqnarray}\label{eqL3.1}
\sum_{n=1}^{\infty} n^{r-2} (\log n)^a \pr\Big(\Delta_n\ge
\sqrt{n}/(\log n)^2,
 M_n^{\ast}\le \frac{\phi(n)}{\sqrt{r-1-\delta}}\Big)
\le C \ep X^2,
\end{eqnarray}
whenever $\delta>0$ is small enough. Let $\beta_n=n
\ep[|X|I\{|X|>\sqrt{n}/n^p\}]$. Then
$|\ep\sum_{i=1}^jX_{ni}^{\prime}|\le \beta_n$, $1\le j\le n$.
Setting
$$\Cal L=\{n:\beta_n\le \frac 18 \sqrt{n}/(\log  n)^2\},$$
we have
$$\{\Delta_n\ge \sqrt{n}/(\log  n)^2\}
\subset \bigcup_{j=1}^n \{ X_j\ne X_{nj}^{\prime}\}, \quad n\in
\Cal L. $$ So for $n\in \Cal L$,
\begin{eqnarray*}
I_n^{\prime}: &=&\pr\Big(\Delta_n\ge \sqrt{n}/(\log  n)^2,
 M_n^{\ast}\le \frac{\phi(n)}{\sqrt{r-1-\delta}}\Big)\\
&\le& \sum_{j=1}^n \pr\Big(X_j\ne X_{nj}^{\prime}, M_n^{\ast}\le
\frac{\phi(n)}{\sqrt{r-1-\delta}}\Big).
\end{eqnarray*}
Observer that $X_{nj}^{\prime}=0$ whenever $X_j\ne
X_{nj}^{\prime}$, $j\le n$, so that by Lemma \ref{lem1}(a) we have
for $n$ large enough and all $1\le j\le n$,
\begin{eqnarray*}
& &\pr\Big(X_j\ne X_{nj}^{\prime},
M_n^{\ast}\le \frac{\phi(n)}{\sqrt{r-1-\delta}}\Big)\\
&=&\pr\Big(X_j\ne X_{nj}^{\prime}, \max_{k\le
j-1}|S_{nk}^{\ast}|\vee \max_{j<k\le
n}|S_{nk}^{\ast}-X_{nj}^{\prime}|
\le \frac{\phi(n)}{\sqrt{r-1-\delta}}\Big)\\
&=&\pr(X_j\ne X_{nj}^{\prime}) \pr\Big(\max_{k\le
j-1}|S_{nk}^{\ast}|\vee \max_{j<k\le
n}|S_{nk}^{\ast}-X_{nj}^{\prime}|
\le \frac{\phi(n)}{\sqrt{r-1-\delta}}\Big)\\
&\le &\pr\Big(X_j\ne X_{nj}^{\prime}) \pr\Big(M_n^{\ast}\le
\frac{\phi(n)}{\sqrt{r-1-\delta}}
+|X_{nj}^{\prime}|\Big)\\
&\le &\pr(|X|>\sqrt{n}/n^p) \pr\Big(M_n^{\ast}\le
\frac{\phi(n)}{\sqrt{r-1-\delta}}
+\sqrt{n}/n^p\Big)\\
&\le &C\pr(|X|>\sqrt{n}/\log^p
n)\exp\{(-(r-1)+\delta+\delta^{\prime})\log  n\}\\
 &\le
&C\ep[|X|^{2+\epsilon}]
n^{-(r-1)+\delta+\delta^{\prime}-(1/2-p)(2+\epsilon)}.
\end{eqnarray*}
Notice that $(1/2-p)(2+\epsilon)\ge 1+\epsilon/4$. Choose
$0<\delta, \delta^{\prime}< \epsilon/16$. Then
$\delta+\delta^{\prime}-(1/2-p)(2+\epsilon)<-1-\epsilon/8$. So,
\begin{eqnarray*}
 \sum_{n\in \Cal L} n^{r-2} (\log n)^a I_n^{\prime} \le
C\ep[|X|^{2+\epsilon}]\sum_{n=1}^{\infty}n^{\delta+\delta^{\prime}-(1/2-p)(2+\epsilon)}(\log
n)^a \le  C\ep[|X|^{2+\epsilon}].
\end{eqnarray*}
Note that
$$ \frac{\beta_n}{\sqrt{n}/(\log n)^2}\le C
\ep[|X|^{2+\epsilon}](\log n)^2 n^{1/2-(1/2-p)(1+\epsilon)}\to 0.
$$
It follows that there are only finite many $n$s not in $\Cal L$.
So
\begin{eqnarray*}
 \sum_{n\not\in \Cal L} n^{r-2}  (\log n)^aI_n^{\prime}<\infty.
\end{eqnarray*}
(\ref{eqL3.1}) is proved.

\begin{lemma}\label{lem4}
Suppose $\ep[|X|^{2+\epsilon}]<\infty$ for some $0<\epsilon<1$.
Let $a>-1$, $r>1$ and $0<p<\epsilon/(4(2+\epsilon))$. Then there
exist constants $\delta_0=\delta_0(r,a,p,\epsilon)>0$ and
$K=K(r,a,p,\epsilon,\delta_0)$ such that $\forall
0<\delta<\delta_0$,
$$\sum_{n=1}^{\infty} n^{r-2} (\log n)^a II_n\le K <\infty, $$
where
$$II_n=\pr\Big(\Delta_n\ge \sqrt{n}/(\log  n)^2,
 M_n\le \phi(n)\big(\frac 1{\sqrt{r-1}}+\delta\big)\Big). $$
\end{lemma}
{\bf Proof.} It is enough to show that
\begin{eqnarray*}
\sum_{n=1}^{\infty}  n^{r-2} (\log n)^a  \pr\Big(\Delta_n\ge
\sqrt{n}/(\log n)^2,
 M_n\le \frac{\phi(n)}{\sqrt{r-1-\delta}}\Big)
\le C \ep |X|^{2+\epsilon},
\end{eqnarray*}
whenever $\delta>0$ is small enough. Let $\beta_n$ and $\Cal L$ be
defined as in the proof of Lemma \ref{lem3}. Then for $n\in \Cal
L$,
\begin{eqnarray*}
\pr\Big(\Delta_n\ge \sqrt{n}/(\log  n)^2,
 M_n\le \frac {\phi(n)}{\sqrt{r-1-\delta}}\Big)
\le \sum_{j=1}^n \pr\Big(X_j\ne X_{nj}^{\prime}, M_n\le
\frac{\phi(n)}{\sqrt{r-1-\delta}}\Big).
\end{eqnarray*}
and for $1\le j\le n$,
\begin{eqnarray*}
& &\pr\Big(X_j\ne X_{nj}^{\prime}, M_n\le \frac{\phi(n)}{\sqrt{r-1-\delta}}\Big)\\
&\le&\pr\Big(X_j\ne X_{nj}^{\prime},
M_{j-1}\vee \max_{j<k\le n}|S_k-X_j+X_j|\le \frac{\phi(n)}{\sqrt{r-1-\delta}},M_n\le \frac{\phi(n)}{\sqrt{1+a-\delta}}\Big) \\
&\le&\pr\Big(\frac{\sqrt{n}}{n^p}<|X_j| \le
2\frac{\phi(n)}{\sqrt{r-1-\delta}},
M_{j-1}\vee \max_{j<k\le n}|S_k-X_j+X_j|\le \frac{\phi(n)}{\sqrt{r-1-\delta}}\Big) \\
&=&\int_{\frac{\sqrt{n}}{n^p}<|y|\le
2\frac{\phi(n)}{\sqrt{r-1-\delta}}} \pr\Big(M_{j-1}\vee
\max_{j<k\le n}|S_k-X_j+y|\le
\frac{\phi(n)}{\sqrt{r-1-\delta}}\Big) d \pr(X_j<y).
\end{eqnarray*}
Note that $M_{j-1}\vee \max_{j<k\le n}|S_k-X_j+y|\overset{\Cal D}=
M_{j-1}\vee \max_{j\le k\le n-1}|S_k+y|$. By Lemma \ref{lem2}, we
have
\begin{eqnarray*}
\sup_{|y|\le 2\frac{\phi(n)}{\sqrt{r-1-\delta}}}
 &\pr&\Big(M_{j-1}\vee \max_{j<k\le n}|S_k-X_j+y|
 \le \frac{\phi(n)}{\sqrt{r-1-\delta}}\Big)\\
&\le& C\exp\big\{(-(r-1)+\delta+\delta^{\prime})\log  n\big\}.
\end{eqnarray*}
It follows that for $n\in \Cal L$ and $1\le j\le n$,
\begin{eqnarray*}
 \pr\Big(X_j\ne X_{nj}^{\prime}, M_n\le \frac{\phi(n)}{\sqrt{r-1-\delta}}\Big)
&\le& C n^{-(r-1)+\delta+\delta^{\prime}}
 \pr\Big(\frac{\sqrt{n}}{n^p}<|X_j|
\le 2\frac{\phi(n)}{\sqrt{1+a-\delta}}\Big)\\
&\le& Cn^{-(r-1)+\delta+\delta^{\prime}}
 \pr\Big(|X|>\sqrt{n}/ n^p\Big).
\end{eqnarray*}
The remained proof is similar to that of (\ref{eqL3.1}) with Lemma
\ref{lem1}(b) instead of Lemma \ref{lem1}(a).

\begin{lemma}\label{lem5}
For any sequence of independent random variables $\{\xi_n; n\ge
1\}$ with mean zero and finite variance, there exists a sequence
of independent normal variables $\{\eta_n; n\ge 1\}$ with $\ep
\eta_n=0$ and $\ep \eta_n^2=\ep \xi_n^2$ such that, for all $Q>2$
and $y>0$,
$$ \pr\Big(\max_{k\le n}|\sum_{i=1}^k \xi_i-\sum_{i=1}^k \eta_i|\ge y\Big)
\le (AQ)^Qy^{-Q}\sum_{i=1}^n \ep |\xi_i|^Q, $$ whenever
$\ep|\xi_i|^Q<\infty$, $i=1,\ldots,n$. Here, $A$ is a universal
constant.
\end{lemma}

{\bf Proof.} See Sakhaneko (1980,1984, 1985).

\begin{lemma}\label{lem6}
We have that
\begin{eqnarray}\label{eqL6.1}
& &\pr\big(\sup_{0\le s\le 1}|W(s)|\le x-1/(\log
n)^2\big)-p_n\le\pr\big(M_n^{\ast}\le x\sqrt{B_n}\big) \non &\le&
 \pr\big(\sup_{0\le s\le 1}|W(s)|\le x+1/(\log  n)^2\big)+p_n,
\quad \forall x>0,
\end{eqnarray}
 where $p_n\ge 0$ satisfies
\begin{equation}\label{eqL6.2}
\sum_{n=1}^{\infty} n^{r-2}(\log n)^a p_n \le K(r,a,p)<\infty.
\end{equation}
\end{lemma}

{\bf Proof.} By Lemma \ref{lem5}, there exist a universal constant
$A>0$ and a sequence of standard Wiener processes $\{W_n(\cdot)\}$
such that for all $Q>2$,
\begin{eqnarray*}
& & \pr\Big( \max_{k\le n}|S_{nk}^{\ast}-W_n(\frac kn B_n)|
\ge \frac 12 \sqrt{B_n}/(\log  n)^2\Big) \\
&\le& (AQ)^Q\Big(\frac{(\log  n)^2}{\sqrt{B_n}}\Big)^Q
\sum_{k=1}^n \ep\big|X_{nk}^{\ast}\big|^Q
\le  C n
\Big(\frac{(\log  n)^2}{\sqrt{n}}\Big)^Q
 \ep\big[|X|^QI\{|X|\le \sqrt{n}/n^p\}\big].
\end{eqnarray*}
On the other hand, by Lemma 1.1.1 of Cs\"org\H o and R\'ev\'esz
(1981),
\begin{eqnarray*}
& & \pr\Big(|\max_{0\le s\le 1}|W_n(s B_n)-W_n(\frac{[ns]}{n}
B_n)|
\ge \frac 12 \sqrt{B_n}/(\log  n)^2\Big) \\
&=& \pr\Big(|\max_{0\le s\le 1}|W_n(s)-W_n(\frac{[ns]}{n} )|
\ge \frac 12 \sqrt{\frac 1n}\frac{\sqrt n}{(\log  n)^2}\Big) \\
&\le& Cn\exp\Big\{-\frac{(\frac 12 \sqrt{n}/(\log  n)^2
)^2}{3}\Big\} \le C n\exp\Big\{-\frac 1{12}n/(\log  n)^4\Big\}.
\end{eqnarray*}
Let
\begin{equation} \label{eqL6.3}
p_n= \pr\Big(\big|M_n^{\ast}/\sqrt{B_n}-\sup_{0\le s\le
1}|W_n(sB_n)|/\sqrt{B_n}\big|\ge 1/(\log  n)^2\Big).
\end{equation}
Then $p_n$ satisfies (\ref{eqL6.1}), since $\{W_n(t
B_n)/\sqrt{B_n}; t\ge 0\}\overset{\Cal D}=\{W(t); t\ge 0\}$ for
each $n$. And also,
$$p_n\le C n \Big(\frac{(\log  n)^2}{\sqrt{n}}\Big)^Q
 \ep\big[|X|^QI\{|X|\le \sqrt{n}/n^p\}\big]
+C n\exp\Big\{-\frac 1{12}n/(\log  n)^4\Big\}. $$
It follows that
\begin{eqnarray*}
&&\sum_{n=1}^{\infty}n^{r-2}(\log n)^a p_n \le K_1+
C\sum_{n=1}^{\infty}n^{r-1-Q/2}(\log n)^{2Q+a}
\ep\big[|X|^QI\{|X|\le \sqrt{n}/n^p\}\big] \\
&\le& K_1+ C\sum_{n=1}^{\infty}n^{r-1-Q/2+(1/2-p)Q}(\log n)^{2Q+a}
 \le K_1+ C\sum_{n=1}^{\infty}n^{r-1-pQ}(\log n)^{2Q+a}
 \le
K<\infty,
\end{eqnarray*}
whenever $Q$ is large enough such that $pQ-r>0$. So,
(\ref{eqL6.2}) is satisfied.

\bigskip
Now, we turn to prove Proposition \ref{prop3.1}. Observe that,
$0\le 1-B_n/n=o((\log n)^{-2})$. If $n$ is large enough, then
\begin{eqnarray*}
& & \pr\Big\{M_n \le \epsilon \phi(n)
 \Big\}\\
&=&\pr\Big\{M_n \le \epsilon \phi(n) , \Delta_n\le
\frac{\sqrt{n}}{(\log  n)^2} \Big\}    +\pr\Big\{M_n \le \epsilon
\phi(n),
\Delta_n> \frac{\sqrt{n}}{(\log  n)^2} \Big\}\\
&\le&\pr\Big\{M_n^{\ast} \le \epsilon \sqrt{\frac{\pi^2 B_n}{8\log
n} } + 2\frac{\sqrt{n}}{(\log  n)^2} \Big\}    +\pr\Big\{M_n \le
\phi(n)\big(\frac 1{\sqrt{r-1}}+\delta\big),
\Delta_n> \frac{\sqrt{n}}{(\log  n)^2} \Big\}\\
&\le&\pr\Big\{\sup_{0\le s\le 1}|W(s)| \le \epsilon
\sqrt{\frac{\pi^2 }{8\log  n} } + \frac{5}{(\log  n)^2} \Big\}
+p_n +II_n
\end{eqnarray*}
for all $\epsilon\in(1/\sqrt{r-1}-\delta,1/\sqrt{r-1}+\delta)$,
where $II_n$ and $p_n$ are defined in Lemmas \ref{lem4} and
\ref{lem6}, respectively. Similarly, if  $n$ is large enough, then
\begin{eqnarray*}
& & \pr\Big\{M_n \le \epsilon \phi(n)\Big\} \ge\pr\Big\{M_n \le
\epsilon \phi(n), \Delta_n\le \frac{\sqrt{n}}{(\log  n)^2} \Big\} \\
&\ge&\pr\Big\{M_n^{\ast} \le \epsilon \phi(n)-
\frac{\sqrt{n}}{(\log  n)^2},
\Delta_n\le \frac{\sqrt{n}}{(\log  n)^2}  \Big\} \\
&\ge&\pr\Big\{M_n^{\ast} \le \sqrt{B_n} \big[\epsilon
\sqrt{\frac{\pi^2 }{8\log  n} } - \frac{2}{(\log  n)^2}\big]
\Big\}
\\
& &\quad -\pr\Big\{M_n^{\ast} \le \epsilon \phi(n)-
\frac{\sqrt{n}}{(\log  n)^2} - \frac{2}{(\log  n)^2}\big],
 \Delta_n> \frac{\sqrt{n}}{(\log  n)^2}\Big\} \\
&\ge&\pr\Big\{\sup_{0\le s\le 1}|W(s)| \le \epsilon
\sqrt{\frac{\pi^2 }{8\log  n} }
 - \frac{3}{(\log  n)^2} \Big\} \\
& & \quad  -\pr\Big\{M_n^{\ast} \le \phi(n)\big(\frac
1{\sqrt{r-1}}+\delta\big),
\Delta_n> \frac{\sqrt{n}}{(\log  n)^2} \Big\}\\
&\ge&\pr\Big\{\sup_{0\le s\le 1}|W(s)| \le \epsilon
\sqrt{\frac{\pi^2 }{8\log  n} } - \frac{5}{(\log  n)^2} \Big\}
-p_n -I_n
\end{eqnarray*}
for all $\epsilon\in(1/\sqrt{r-1}-\delta,1/\sqrt{r-1}+\delta)$,
where $I_n$ is  defined in Lemma \ref{lem3}. Choosing $\delta>0$
small enough and letting $q_n=p_n+I_n+II_n$ complete the proof by
Lemmas \ref{lem3}, \ref{lem4} and \ref{lem6}.


\section{Proof of the Theorems.}
\setcounter{equation}{0} {\bf Proof of Theorem \ref{th1}:} Suppose
(\ref{eqT1.0}) hold. Without losing of generality, we assume that
$\ep X=0$ and $\ep X^2=1$. Let $\delta>0$, $p>0$ small enough and
$\{q_n\}$ be such that (\ref{eqprop3.1.1}) and (\ref{eqprop3.1.2})
hold. Then
\begin{eqnarray*}
 \lim_{\epsilon\nearrow 1/\sqrt{r-1}}
[\epsilon^{-2}-(r-1)]^{a+1} \sum_{n=1}^{\infty} n^{r-2} (\log n)^a
q_n=0,
\end{eqnarray*}
by (\ref{eqprop3.1.2}). Notice that $a_n(\epsilon)\to 0$. By
(\ref{eqprop3.1.1}), we have that for $n$ large enough,
\begin{eqnarray*}
&&\pr\Big\{\sup_{0\le s\le 1}|W(s)| \le \sqrt{\pi^2/(8\log  n)}
(\epsilon+a_n(\epsilon))
 + \frac{5}{(\log  n)^2 } \Big\}
-q_n
\\
&\le& \pr\Big\{M_n \le \phi(n) (\epsilon+a_n(\epsilon)) \Big\}
\\
&\le&\pr\Big\{\sup_{0\le s\le 1}|W(s)| \le \sqrt{\pi^2/(8\log  n)}
(\epsilon+a_n(\epsilon)) - \frac{5}{ (\log  n)^2 }
\Big\} +q_n, \\
& &  \qquad \forall \epsilon\in (1/\sqrt{r-1}-\delta/2,
1/\sqrt{r-1}+\delta/2).
\end{eqnarray*}
On the other hand, by Proposition \ref{prop2.1},
\begin{eqnarray*}
& & \lim_{\epsilon\searrow \sqrt{r-1}} [\epsilon^{-2}-(r-1)]^{a+1}
\sum_{n=1}^{\infty}n^{r-2} (\log n)^a   \pr \Big\{\sup_{0\le s\le
1}|W(s)|
 \le \sqrt{\pi^2/(8\log n)} (\epsilon
+a_n(\epsilon)) \pm \frac{5}{(\log n)^2 }\Big\}
\\
& & \qquad \qquad \qquad=\frac 2{\pi}\exp\{2\tau (r-1)^{3/2}\}
\Gamma(a+1).
\end{eqnarray*}
(\ref{eqT1.1}) is  now proved.

Conversely, suppose (\ref{eqT1.1}) hold. From Esseen (1968) (see
also Petrov 1995, Page 74 (2.70) ) it is easy to see that for all
$m\ge 1$,
$$\pr(|S_m|\le 2\sqrt{m})\le K\big(\int_{-2\sqrt{m}}^{2\sqrt{m}}d
\widetilde{F}(x)\big)^{-1/2},
$$
where $\widetilde{F}(x)$ is the distribution function of the
symmetrized $X$, and $K$ is an absolute constant. So, if $\ep
X^2=\infty$, then for any $M>2$ we can choose $m_0\ge 9$ large
enough such that
\begin{equation}\label{eq3.1} \pr(|S_m|\le
2\sqrt{m})\le e^{-2M}, \quad m\ge m_0. \end{equation}
 For
$\epsilon>0$, we let $m=[\epsilon^2n/\log n]$, and $N=[n/m]$, then
for all $n\ge m_0^2$ and $\epsilon\ge 1$,
\begin{eqnarray*}
\pr\big(M_n\le \epsilon (n/\log n)^{1/2}\big)
&\le&\pr\big(|S_{km}-S_{(k-1)m}|\le 2\sqrt{m}, k=1,\ldots, N\big)
\non &\le& e^{-2MN}\le \exp\big\{-M\frac{\log
n}{\epsilon^2}\big\}.
\end{eqnarray*}
By this inequality, for any $r$, $a>-1$ and
$0<\epsilon_1<\epsilon_2<\infty$ there exists a constant
$C=C(r,a,\epsilon_1,\epsilon_2)$ for which
$$ \sup_{\epsilon\in (\epsilon_1,\epsilon_2)}
\sum_{n=1}^{\infty}n^{r-2}(\log n)^a
 \pr\big\{M_n\le \epsilon (n/\log n)^{1/2}\big\}\le C<\infty,$$
which implies that
$$ \lim_{\epsilon\nearrow 1/\sqrt{r-1}}[\epsilon^{-2}-(r-1)]^{a+1}
\sum_{n=1}^{\infty}n^{r-2}(\log n)^a
 \pr\big\{M_n\le \phi(n)(\epsilon+a_n(\epsilon))\big\}=0.$$
This contradictory to (\ref{eqT1.1}). If $\ep X^2<\infty$ and $\ep
X=\mu\ne 0$, then (\ref{eq3.1}) also hold since $|S_m|/\sqrt{m}\to
\infty$ a.s. as $m\to \infty$. So, we conclude   $\ep X^2<\infty$
and $\ep X=0$. At last, under $\ep X^2<\infty$ and $\ep X= 0$,
$\Var X=\sigma^2$ is obvious according (\ref{eqT1.1}) and Theorem
\ref{th0}.

\bigskip
{\bf Proof of Theorem \ref{th3}:} Let $\delta>0$, $p>0$ small
enough and $\{q_n\}$ be such that (\ref{eqprop3.1.1}) and
(\ref{eqprop3.1.2}) hold. By a standard argument (see Feller
(1945)), we can assume that
$$ \big(\frac{1}{\sqrt{r-1}}-\frac{\delta}{2}\big)\phi(n)
\le \sqrt{\frac{\pi^2 n}{8\psi(n)}}\le
\big(\frac{1}{\sqrt{r-1}}+\frac{\delta}{2}\big)\phi(n).
$$
That is
$$
\big(\frac{1}{\sqrt{r-1}}+\frac{\delta}{2}\big)^{-1/2}\log n \le
\psi(n)\le
\big(\frac{1}{\sqrt{r-1}}-\frac{\delta}{2}\big)^{-1/2}\log n.
$$
Let $\epsilon=\sqrt{\pi^2n/(8\psi(n)}/\phi(n)$. By
(\ref{eqprop3.1.1}), it follows that
\begin{eqnarray*}
&&\pr\Big\{\sup_{0\le s\le 1}|W(s)| \le \sqrt{\pi^2/(8\psi(n))}
 + 5/(\log  n)^2  \Big\}
-q_n
\\
&\le& \pr\Big\{M_n \le \sqrt{\pi^2 n /(8\psi(n))} \Big\}
\\
&\le&\pr\Big\{\sup_{0\le s\le 1}|W(s)| \le \sqrt{\pi^2/(8\psi(n))}
 - 5/(\log  n)^2  \Big\}
+q_n.
\end{eqnarray*}
Notice
$$\frac{\pi^2}{ 8\big(\sqrt{\pi^2/(8\psi(n))} \pm 5/
(\log n)^2 \big)^2}=\psi(n)\big(1+o(\frac{1}{\log
n})\big)=\psi(n)+o(1).
$$
According to (\ref{eqL2.1.1}), it follows that
$$ c_1\exp\{-\psi(n)\}-q_n\le
\pr\Big\{M_n \le \sqrt{ \pi^2 n /(8\psi(n))} \Big\} \le c_2
\exp\{-\psi(n)\}+q_n.
$$
By (\ref{eqprop3.1.2}), Theorem \ref{th3} is now proved.

\newpage


\end{document}